\documentclass[12pt]{article}
\usepackage{latexsym}
\usepackage{latexsym}
\usepackage[labelfont=bf]{caption}
\newtheorem{thm}{Theorem}[section]

\newtheorem{cor}[thm]{Corollary}
\newtheorem{prop}[thm]{Proposition}
\newtheorem{defn}{Definition}[section]

\newtheorem{rem}{Remark}[section]

\newfont{\sss}{cmss12 scaled 1000}

\title{Construction of 3-Designs Using $(1,\sigma)$-Resolution}
\author{ Tran van Trung \\
    Institut f\"ur Experimentelle Mathematik \\
    Fakult\"at f\"ur Mathematik \\ 
     Universit\"at Duisburg-Essen \\
    Thea-Leymann-Stra\ss e 9, 45127 Essen, Germany }

\date{}

\begin{document}

\maketitle

\begin{abstract}
 The paper deals with recursive constructions for simple 3-designs based 
 on other 3-designs having $(1, \sigma)$-resolution. The concept of 
 $(1, \sigma)$-resolution may be viewed as a generalization of the 
 parallelism for designs. We show the constructions and their applications
 to produce many previously unknown infinite families of simple 3-designs.
 We also include a discussion of 
 $(1,\sigma)$-resolvability of the constructed designs.
\end{abstract}

\vspace{0.1in}\noindent

\vspace{0.1in}\noindent

\bigskip
\section{Introduction}

 In our previous papers \cite{tvt2001, tvt2000} we have presented
 several recursive constructions for simple 3-designs.
 In \cite{tvt2001}, among others, generalizations
 of the well-known doubling construction of Steiner quadruple
 systems for 3-designs  are introduced. In \cite{tvt2000} more general
 recursive constructions of simple 3-designs are described, 
 whereby ingredient designs may have repeated blocks.
 The methods in these papers are based on the existence of 
 3-designs having a parallelism, i.e. the blocks of the design
 can be partitioned into classes of mutually disjoint blocks 
 such that every point is in exactly one block
 of each class. Designs with parallelism have shown to be 
 useful for constructing designs in the literature
 \cite{Shrikhande63}, \cite{driessen78}, \cite{jv86}, 
 \cite{Phelps89}, \cite{Jimbo2011}, \cite{Laue2004}, \cite{Stinson2014},
 \cite{tvt2001, tvt2000}.
  
 The concept of $(1,\sigma)$-resolvability for $t-(v,k,\lambda)$ 
 designs may be viewed as a generalization of that of parallelism.  
 For the latter means that the design is $(1,1)$-resolvable. It should
 be mentioned that if a $t-(v,k,\lambda)$ design has a parallelism we 
 necessarily have $k | v$; this condition does no longer hold 
 for $(1,\sigma)$-resolvability in general.  
 Thus, the natural question is that whether or not the 
 methods in our previous papers \cite{tvt2001, tvt2000} can be extended to 
 $(1,\sigma)$-resolvable 3-designs. We show that this is in fact the case.
 Our aim in this paper is to present this generalization. 
 The result provides a general method for constructing simple 3-designs
 which largely extends the use of complete designs as ingredients
 for the construction. We show the strength of the method by  
 giving some simple applications to construct a number of families of 
 simple 3-designs, which, to our knowledge, were not
 previously known to exist. We also include a discussion
 of $(1, \sigma)$-resolvability of the constructed designs.  
 
 For notation and general definitions of $t$-designs we refer to 
 \cite{bjl86, handbook96}.


\bigskip

 \section{Constructions of 3-Designs using $(1,\sigma)$-Resolution}

  In this section we present recursive constructions of simple
  3-designs using $(1,\sigma)$-resolution of their ingredients. 
 
 \subsection{Preliminaries} 

  \medskip 
  We begin with a few definitions and set up necessary conditions 
  for the ingredients used in the constructions.   

 \begin{defn}\label{def1}
  A $t-(v,k,\lambda)$-design $(X, {\mathcal B})$ is
 said to be $(s,\sigma)$-resolvable for a given $s \in \{ 1, \ldots, t \}$, 
 if its block set $\mathcal B$ 
 can be partitioned into $w$ classes $\pi_1, \ldots, \pi_w$ 
 such that
 $(X,\pi_i)$ is a $s-(v,k,\sigma)$ design for all $i=1,\ldots, w.$
 Each $\pi_i$ is called a resolution class.
 \end{defn}
 It is worth noting that the concept of resolvability 
 (i.e. $(1,1)$-resolvability) for BIBD introduced
 by Bose in 1942 \cite{Bose1942} was generalized by
 Shrikhande and Raghavarao to $\sigma$-resolvability (i.e. 
 $(1, \sigma)$-resolvability) for BIBD in 1964 \cite{SR1963}.
 A definition of $s$-resolvability (i.e. $(s, \sigma)$-resolvability) 
 for t-designs with $t \geq 3$ and $1 \leq s \leq t$ may be found in
 \cite{Baker1976}, for example.
 
 \begin{rem}
   If $(X, {\mathcal B})$ is the complete $t-(v,k, {v-t\choose k-t})$
  design, then a $(t, \sigma)$-resolution of $(X, {\mathcal B})$ is
  a large set of $t-(v,k,\sigma)$ designs. 
  
  It should be remarked that
  each $t-(v,k,\lambda)$ design always has a {\it trivial} 
  $(s,\lambda_s)$-resolution 
  consisting of a single class, i.e. $w=1$, for all
  $1 \leq s \leq t$. Throughout the paper when we speak of
  $(s, \sigma)$-resolution we mean that $w \geq 2.$
  Note that  
  $w=\lambda{v \choose t} {k\choose s}/\sigma{v \choose s}{k\choose t}$.
 \end{rem}

 \begin{defn}
  Let $D$ be a $t-(v,k,\lambda)$ design admitting
  a $(s,\sigma)$-resolution with
  $\pi_1, \ldots , \pi_w$ as resolution classes.
  Define a distance  between any two classes 
  $\pi_i$ and $\pi_j$ by
  $d(\pi_i, \pi_j)= \min \{|i-j|, w-|i-j| \}.$    
 \end{defn}
 
  \medskip
  For the constructions in this paper we employ designs 
  having a $(1,\sigma)$-resolution. We now describe the
  detailed assumption and notation used throughout the paper.

  \medskip 
 
  Let $\{k_1, \ldots, k_n, k_{n+1}, \ldots, k_{2n} \}$ and 
  $k$ be integers
  with $2\leq k_1 < \cdots < k_n  \leq k/2$  such that
  $k_i+ k_{n+i}=k$ for $i=1, \ldots, n$.
 
  Assume that there exist $3-(v,k_i, \lambda^{(i)})$
  designs $D_i=(X,{\mathcal B}_i)$ having a $(1,\sigma^{(i)})$-resolution
  such that $w_i= w_{n+i}$ for all $i=1, \ldots, n$, 
  where $w_j$ denotes the number of 
  classes in a $(1,\sigma^{(j)})$-resolution of $D_j$, i.e.  
  $D_i$ and $D_{n+i}$ have the same number of resolution classes.  
  
  It is also assumed that 
 \begin{enumerate} 
 \item 
  For each pair $(D_i, D_{n+i})$, $1\leq i \leq n$,
  either $D_i$ or $D_{n+i}$ has to be simple.
 \item
  If a $D_j$, $j\in \{i,n+i\}$, is not simple, then $D_j$
  is a union of $a_j$ copies of a simple
  $3-(v,k_j, \alpha^{(j)})$ design $C_j$, wherein $C_j$ admits
  a $(1,\sigma^{(j)})$-resolution. Thus, 
  $\lambda^{(j)}= a_j \alpha^{(j)}$. 
 \end{enumerate}

  Note that the trivial $2-(v,2,1)$ design will be 
  considered as a $3-(v,2,\lambda)$ design with $\lambda=0$.  
 
 \medskip
  Further we need to specify the way of setting up 
 $(1,\sigma^{(j)})$-resolution
  classes for $D_j$, when $D_j$ is the union of $a_j$ copies $C_j$. 
  
  Let $P^{(j)}= \{ \pi^{(j)}_1, \ldots , \pi^{(j)}_{t_j} \}$ 
  be a $(1,\sigma^{(j)})$-resolution
  of the simple design $C_j$. The corresponding
  $(1,\sigma^{(j)})$-resolution of $D_j$ is
  chosen to be the ``concatenation'' of $a_j$ sets
  $P^{(j)}$. This means that the $w_j=a_j t_j$ resolution classes 
  of $D_j$ are arranged in the following way
  \[ \pi^{(j)}_1, \ldots , \pi^{(j)}_{t_j},\ \ 
     \pi^{(j)}_1, \ldots , \pi^{(j)}_{t_j}, \ 
    \ldots \ , \ \ \pi^{(j)}_1, \ldots, \pi^{(j)}_{t_j}. \]
              
  Finally, we also assume that there exists
  a $3-(v,k,\Lambda)$ design $D=(X,{\mathcal B})$, when it is 
  needed in our construction. 
  
  \medskip
  \noindent
   {\bf Notation: \quad}
   \begin{itemize}
   \item     
    $\pi^{(\ell)}_1, \ldots , \pi^{(\ell)}_{w_{\ell}}$ 
  denote the $w_{\ell}$ classes in a $(1,\sigma^{(\ell)})$-resolution
  of $D_{\ell}$
  for $\ell=1,\ldots, 2n$. Recall that $w_h= w_{n+h}$ for $h=1,\ldots,n.$
    \item
   The distance defined on the resolution classes of $D_{\ell}$ 
   is then $d^{(\ell)}(\pi^{(\ell)}_i, \pi^{(\ell)}_j)= 
    \min \{|i-j|, w_{\ell}-|i-j| \}.$   

   \item
    $b^{(j)}=\sigma^{(j)}v/k$  denotes
    the number of blocks in each class of 
    a $(1,\sigma^{(j)})$-resolution of $D_j$.

   \item 
    $u_j:=\sigma^{(j)}$  denotes
    the number of blocks containing a point in each class of 
    a $(1,\sigma^{(j)})$-resolution of $D_j$.

    \item
    $\lambda^{(j)}_2=\lambda^{(j)}(v-2)/(k_j-2)$  denotes
    the number of blocks of $D_j$ containing two points.

 \end{itemize}

\subsection{Construction I}
 
 \medskip 
 In this section we describe the first construction by using the 
 set-up above for the case $k_n \not= k/2.$ 

\medskip
 
  Let $\tilde D_i=(\tilde X,\tilde{\mathcal B}_i)$ be a copy of $D_i$ 
  defined on the point set $\tilde X$ such that $X \cap \tilde X=\emptyset$.  
 Also
 let $\tilde D=(\tilde X,\tilde{\mathcal B})$ be a copy of $D$.
  
   Define blocks on the point set $X\cup \tilde X$ as follows:
 \begin{itemize}
 
 \item[I.] blocks of $D$ and blocks of  $\tilde D$;  
 
 \item[II.] blocks of the form $A \cup \tilde B$ for any 
     $A\in \pi^{(h)}_i$ and $\tilde B \in {\tilde \pi}^{(n+h)}_j$ with
  $\varepsilon_h \leq d^{(h)}(\pi^{(h)}_i, \pi^{(h)}_j)\leq s_h$, 
  $\varepsilon_h =0,1$,
  for $h=1,\ldots, n$;
  
 \item[III.] blocks of the form $\tilde A \cup B$ for any 
     $\tilde A\in {\tilde \pi}^{(h)}_i$ and $B \in \pi^{(n+h)}_j$ with
  $\varepsilon_h \leq d^{(h)}(\pi^{(h)}_i, \pi^{(h)}_j)\leq s_h$, 
  $\varepsilon_h =0,1$, for $h=1,\ldots, n$.
   
 \end{itemize}
 
 \medskip
 Here, and in the sequel, 
 the non-negative integers $s_h$, $h=1,\ldots, n$, denote the parameters 
 that have to be determined, for which the
 defined blocks of types I, II and III form a 3-design. Thus, $s_h$,
 should not be confused with $s$ in $(s, \sigma)$-resolution as defined
 above.

\bigskip
\noindent
 {\bf  Notation:} Define $z_h=(2s_h+1-\varepsilon_h)$ 
   if $s_h < \frac{w}{2}$, and $z_h=(2s_h-\varepsilon_h)$
   if $s_h = \frac{w}{2}$, for $h=1,\ldots, n$.

 \bigskip 
 Any 3 points $a,b,c\in X$, resp. 
 $\tilde a,\tilde b,\tilde c \in \tilde X$   
 are contained in 
 
 \begin{itemize}
 
  \item
  $\Lambda$ blocks of type I,
 
  \item  
  $z_h\lambda^{(h)}b^{(n+h)}$ blocks 
  of type II for $h=1,\ldots, n$,
  
  \item
   $z_h\lambda^{(n+h)}b^{(h)}$ blocks of type III
   for $h=1,\ldots, n$.
 
 \end{itemize}
   
 Thus $a,b,c $ appear together in
  $$ \Lambda + \sum_{h=1}^n 
    z_h\lambda^{(h)}b^{(n+h)}+ z_h\lambda^{(n+h)}b^{(h)}$$
   blocks. Set
   $$\Delta = 
   \sum_{h=1}^n z_h\lambda^{(h)}b^{(n+h)}+ z_h\lambda^{(n+h)}b^{(h)}.$$
 
\medskip  
 Now consider 3 points $a,b,\tilde c$, where $a,b\in X$
 and $\tilde c \in \tilde X$. Because of the symmetry the number 
 of blocks containing 3 points $a,b,\tilde c$, 
 is equal to the number of blocks containing
 $\tilde a,\tilde b,c$.
 For each $h=1,\ldots, n$, any two points $a$ and $b$
 are contained in $\lambda^{(h)}_2$ blocks of $D_h$ 
 and in $\lambda^{(n+h)}_2$ blocks of $D_{n+h}$; 
 further, the point $\tilde c$ is in  
 $u_h$ (resp. $u_{n+h}$) blocks 
 of each resolution class of $\tilde D_h$ (resp. $\tilde D_{n+h}$).  

 So $a,b,\tilde c$ appear in
 
 \begin{itemize}
 
 \item
 $z_h \lambda^{(h)}_2 u_{n+h}$ blocks of type II
  for $h=1,\ldots, n$,
 \item
 $z_h \lambda^{(n+h)}_2 u_h$ blocks 
 of type III for $h=1,\ldots, n$.
 
 \end{itemize}
 
 Thus $a,b,\tilde c$ are contained together in 
 $$ \Theta:=\sum_{h=1}^n 
 z_h \lambda^{(h)}_2 u_{n+h} + z_h \lambda^{(n+h)}_2 u_h $$
  blocks.

\medskip
 Therefore the blocks defined in I, II and III will form a 3-design if

 \[\Lambda + \Delta = \Theta, \]
or
 \[\Lambda = \Theta - \Delta. \]

 Note that $\Lambda = \Theta - \Delta \geq 0$. The case 
 $\Lambda=\Theta - \Delta = 0$ implies that               
 $D$ and $\tilde D$ are not needed in the construction.
 In both cases either $\Theta - \Delta > 0$ or $\Theta - \Delta =0$ 
 the constructed blocks form a simple
 $3-(2v,k,\Theta)$ design with 
   \[\Theta= 
    \sum_{h=1}^n \{(\lambda^{(h)}_2 u_{n+h} + 
    \lambda^{(n+h)}_2 u_h) \}z_h, \]  
   where  $1\leq z_h \leq w_h$ if both $D_h$ and $D_{n+h}$ are simple
   and $1\leq z_h \leq t_j$ if $D_j$ is non-simple, $j \in \{h, n+h\}.$

 \bigskip
 
   What remains to be verified is the  
   simplicity of the resulting design when 
   either $D_h$ or $D_{n+h}$ is non-simple.
   Evidently, if both $D_h$ and $D_{n+h}$ are simple for all
   $ 1\leq h\leq n$, then the constructed design is simple.
   
 \bigskip
   
 To start with we observe that two blocks
 constructed from two pairs $(D_i, D_{n+i})$ and  $(D_j, D_{n+j})$,
 $i\not= j$, are always distinct. 
 Further any two blocks of different types 
 are also distinct. Thus, we need to consider
 two blocks of the same type, in particular, of type II or type III  
 constructed from a pair $(D_j, D_{n+j})$. 
 W.l.o.g. we may assume that $D_j$ is a union 
 of $a_j$ copies of a simple
 $3-(v,k_j, \alpha^{(j)})$ design $C_j$ and $D_{n+j}$ is simple.

 The following argument is the same for blocks of types II
 and III. So let 
 $ E= A_1\cup \tilde B_1$ and $ F= A_2\cup \tilde B_2$  
 be two blocks of type II of the resulting design, where
 $A_1 \in \pi^{(j)}_{i_1}$, $\tilde B_1 \in \tilde \pi^{(n+j)}_{h_1}$,
 $A_2 \in \pi^{(j)}_{i_2}$ and $\tilde B_2 \in \tilde \pi^{(n+j)}_{h_2}$.
 Suppose $E=F$. Then $\tilde B_1= \tilde B_2$, and hence $h_1=h_2$, since
 $\tilde D_{n+j}$ is simple. Consequently, $A_1= A_2$, so we have 
 
 \begin{enumerate}
 \item either
 $i_1=i_2$,   
 
  \item or
  $i_1\not=i_2$.    
 \end{enumerate}
 
 In the first case, $E$ and $F$ are the same block.
 In the second case, $E$ and $F$ are repeated blocks; this can happen
 only if  $|i_2-i_1|$ is a multiple of $t_j$, i.e.
 $t_j| \ |i_2-i_1|$, this is because the  
 resolution classes of $D_j$ are chosen to be the
 concatenation of $a_j$ copies of a given set $P^{(j)}$ of
 resolution classes of $C_j$. Now, as
 $\varepsilon_j \leq d^{(j)}(\pi^{(j)}_{i_1}, \pi^{(j)}_{h_1}) \leq s_j$
 and
 $\varepsilon_j \leq d^{(j)}(\pi^{(j)}_{i_2}, \pi^{(j)}_{h_2}) 
 =d^{(j)}(\pi^{(j)}_{i_2}, \pi^{(j)}_{h_1}) \leq s_j$,
 it follows that $ z_j > t_j. $ Therefore, the second
 case will not occur if $z_j \leq  t_j.$
 
 Hence, if $z_j \leq t_j$ for
 all non-simple $D_j$'s, the resulting design remains simple.

       
\vspace{2mm}
   With the notation above,    
   we summarize Construction I in the following theorem.

 \begin{thm}\label{Construction1}
  Let $\{k_1, \ldots, k_n, k_{n+1}, \ldots , k_{2n} \}$ 
  and $k$ be integers
  with $2\leq k_1 < \cdots < k_n < k/2$ and $k_i+ k_{n+i}=k$
  for $i=1, \ldots, n$.   
  Assume that there exist  $3-(v,k_i, \lambda^{(i)})$
  designs $D_i=(X,{\mathcal B}_i)$ admitting a 
  $(1,\sigma^{(i)})$-resolution such that $w_i=w_{n+i}$, where 
  $w_j$ is the number of resolution classes of $D_j$.
  Assume further that at least one design from each 
  pair $(D_i, D_{n+i})$, $1\leq i \leq n$, is simple and
  if a $D_j$, $j\in \{i,n+i\}$, is not simple, then $D_j$
  is a union of $a_j$ copies of a simple
  $3-(v,k_j, \alpha^{(j)})$ design $C_j$
  admitting a $(1,\sigma^{(j)})$-resolution, i.e. 
  $\lambda^{(j)}= a_j \alpha^{(j)}$. Let $t_j$ denote the number of
  resolution classes of $C_j$. 
  Let
  $$\Theta := \sum_{h=1}^n \{(\lambda^{(h)}_2 u_{n+h} + \lambda^{(n+h)}_2 u_h)\}
  z_h,$$ 
 $$\Delta := \sum_{h=1}^n \{(\lambda^{(h)}b^{(n+h)}+ \lambda^{(n+h)}b^{(h)} )\}     z_h. $$
  
  \begin{itemize}
  
 \item[(i)]
   Assume that
  \begin{eqnarray}
  0 & = & \Theta - \Delta,
  \end{eqnarray}          
   with $1\leq z_h \leq w_h$ if both $D_h$ and $D_{n+h}$ are simple
   and $1\leq z_h \leq t_j$ if $D_j$ is non-simple, $j \in \{h,n+h\}$.
   Then there exists a simple $3-(2v,k,\Theta)$ design $\mathcal D$.

   \item[(ii)]
   Assume that
    \begin{eqnarray}
    0 & < & \Theta - \Delta, 
    \end{eqnarray}                    
   with $1\leq z_h \leq w_h$ if both $D_h$ and $D_{n+h}$ are simple
   and $1\leq z_h \leq t_j$ if $D_j$ is non-simple, $j \in \{h,n+h\}$;
   further  assume that there is a $3-(v,k,\Lambda)$ design 
   with $\Lambda = \Theta - \Delta.$ 
   Then there exists a simple $3-(2v,k,\Theta)$ design $\mathcal D$. 

 \end{itemize}
 \end{thm} 
  
\bigskip

\subsection{Construction II}
 
 \medskip
 In this section we consider the case $k_n=k/2.$

 We observe that the resulting designs in Construction I 
 would have repeated blocks if $k_n=k/2$ and the block sets of 
 $D_n$ and $D_{2n}$ are not disjoint. 
 To deal with the case  $k_n=k/2$ the blocks constructed from
 the pair $(D_n, D_{2n})$ need to be modified.

\medskip 
 Suppose now $2\leq k_1 < \cdots < k_n = k/2$. 
 Take $D_n=D_{2n}$ and assume that $D_n$ is simple. 
 Now define the blocks on the point set
 $X\cup \tilde X$ as follows:
 
 \begin{itemize}
 
 \item[I.] blocks of $D$ and blocks of  $\tilde D$;  
 
 \item[II.] blocks of the form $A \cup \tilde B$ for any 
     $A\in \pi^{(h)}_i$ and $\tilde B \in {\tilde \pi}^{(n+h)}_j$ with
  $\varepsilon_h \leq d^{(h)}(\pi^{(h)}_i, \pi^{(h)}_j)\leq s_h$, 
  $\varepsilon_h =0,1$,
  for $h=1,\ldots, n-1$;
  
 \item[III.] blocks of the form $\tilde A \cup B$ for any
     $\tilde A\in {\tilde \pi}^{(h)}_i$ and $B \in \pi^{(n+h)}_j$ with
  $\varepsilon_h \leq d^{(h)}(\pi^{(h)}_i, \pi^{(h)}_j)\leq s_h$, 
  $\varepsilon_h =0,1$,
  for $h=1,\ldots, n-1$;
   
  \item[IV.] blocks of the form $A \cup \tilde B$ for any 
     $A\in \pi^{(n)}_i$ and $\tilde B \in {\tilde \pi}^{(2n)}_j$ with
  $\varepsilon_n \leq d^{(n)}(\pi^{(n)}_i, \pi^{(n)}_j)\leq s_n$, 
  $\varepsilon_n =0,1$.

 \end{itemize}
  
  Construction II differs from Construction I only in  blocks of type IV. 
  Observe that
  any three points $a,b,c \in X$ (resp. $ \tilde a, \tilde b, \tilde c
  \in  \tilde X$) are contained in 
  $z_n\lambda^{(n)} b^{(n)}$ blocks
  of type IV; any three points $a,b,\tilde c$ with  
  $ a, b \in X$ and $ \tilde c \in  \tilde X$ 
 (resp. $\tilde a, \tilde b, c$) are contained in 
 $z_n\lambda^{(n)}_2 u_n$ blocks of type IV. 
  All other countings as well as the proof of simplicity 
  of the resulting design remain unchanged as shown in Construction I. 
  
 \medskip 
  We obtain the following theorem for the case $k_n=k/2$.
         

 \begin{thm}\label{Construction2}
   Let $\{k_1, \ldots,k_n, k_{n+1}, \ldots, k_{2n} \}$ and $k$ be integers
  with $2\leq k_1 < \ldots < k_n=k/2$ and 
  $k_i+k_{n+i}=k$ for $i=1, \ldots, n$.
   Assume that there exist  $3-(v,k_i, \lambda^{(i)})$
  designs $D_i=(X,{\mathcal B}_i)$ admitting a 
  $(1,\sigma^{(i)})$-resolution such that $w_i=w_{n+i}$, where 
  $w_j$ is the number of resolution classes of $D_j$.
  Assume further that at least one design from each 
  pair $(D_i, D_{n+i})$, $1\leq i \leq n$, is simple and
  if a $D_j$, $j\in \{i,n+i\}$, is not simple, then $D_j$
  is a union of $a_j$ copies of a simple
  $3-(v,k_j, \alpha^{(j)})$ design $C_j$
  admitting a $(1,\sigma^{(j)})$-resolution, i.e. 
  $\lambda^{(j)}= a_j \alpha^{(j)}$. Let $t_j$ denote the number of
  resolution classes of $C_j$. 
  Let
 $$\Theta^* := \lambda_2^{(n)}u_nz_n + 
              \sum_{h=1}^{n-1}\{(\lambda^{(h)}_2 u_{n+h}+ 
              \lambda^{(n+h)}_2 u_h)\}z_h,$$ 
 $$\Delta^* := \lambda^{(n)}b^{(n)}z_n + 
            \sum_{h=1}^{n-1}\{(\lambda^{(h)}b^{(n+h)}+ 
            \lambda^{(n+h)}b^{(h)} )\}z_h. $$
  
  \begin{itemize}
  
 \item[(i)]
   Assume that
  \begin{eqnarray}
  0 & = & \Theta^* - \Delta^*,
  \end{eqnarray}          
   with $1\leq z_h \leq w_h$ if both $D_h$ and $D_{n+h}$ are simple
   and $1\leq z_h \leq t_j$ if $D_j$ is non-simple, $j \in \{h,n+h\}.$
   Then there exists a simple $3-(2v,k,\Theta^*)$ design $\mathcal D$.

   \item[(ii)]
   Assume that
    \begin{eqnarray}
    0 & < & \Theta^* - \Delta^*, 
    \end{eqnarray}                    
   with $1\leq z_h \leq w_h$ if both $D_h$ and $D_{n+h}$ are simple
   and  $1\leq z_h \leq t_j$ if $D_j$ is non-simple, $j \in \{h,n+h\}$;
   further assume that there is a $3-(v,k,\Lambda)$ design 
   with $\Lambda = \Theta^* - \Delta^*.$ 
   Then there exists a simple $3-(2v,k,\Theta^*)$ design $\mathcal D$.

 \end{itemize}
 \end{thm} 

 \bigskip
 
\section{Applications}

\medskip
 In this section we show applications of Constructions I and II
 for some small values of $n$. It turns out that we can construct
 many new infinite families of simple 3-designs by merely using
 complete designs as ingredients. For these applications we implicitly
 use the following result and observation.

\begin{itemize}
\item {\bf Baranyai's Theorem} \cite{bar75}. 
  The trivial $k-(v,k,1)$ design is (1,1)-resolvable
 (i.e. having a parallelism) if and only if  $k|v.$ 

\medskip
\item {\bf Block orbits}. 
      If $\gcd(v,k)=1$, then the $k-(v,k,1)$ design is
      $(1,k)$-resolvable. The resolution classes are the block
      orbits of a fixed point free automorphism of order $v$.

\end{itemize}

\bigskip
\subsection{Applications of Construction I}

\subsubsection{\boldmath${ n=1}$}

\medskip

 We consider the most simple case of Construction I, namely the case with 
 $n=1$, $k_1=2$ and $k_2=3.$

 \medskip
 Let $ v > 5$ be an integer such that $v \equiv 0 \bmod 2$ and 
 $\gcd(v,3)=1.$ 

 \begin{itemize}

 \item $D_1$ is the union of $a_1=(v-2)/6$ copies of 
  the complete $2-(v,2,1)$ design $C_1$.
  By Baranyai's Theorem  $C_1$ is (1,1)-resolvable,
 and the number of resolution (parallel) classes of $C_1$ is $t_1=(v-1)$. 
 For $D_1$ we have $\lambda^{(1)}=0$,  $\lambda_2^{(1)}=(v-2)/6$,
 $u_1=1$, $b^{(1)}=v/2$  and $w_1=a_1t_1.$  

 \item
  $D_2$ is the complete  $3-(v,3,1)$ design.
 Recall by the observation above that
 $D_2$ admits a $(1,3)$-resolution, which is derived from the block
 orbits of a fixed point-free automorphism of order $v$ on the point set.
 For $D_2$ we have  $\lambda^{(2)}=1$, $\lambda_2^{(2)}=v-2$, $u_2=3$, 
 $b^{(2)}=v$ and $w_2=(v-1)(v-2)/6.$

\item  $D$ is the complete 
       $3-(v,5, \Lambda)= 3-(v,5, {v-3 \choose 2})$ design. 

\end{itemize}

 \medskip
 With the notation of Theorem \ref{Construction1} we can check that
 $$\Lambda = \Theta - \Delta $$
 if $z_1= (v-4)/2,$ where
 $$\Theta =\{ \lambda_2^{(1)}u_2+ \lambda_2^{(2)}u_1\}z_1= 3(v-2)z_1/2,$$
 $$\Delta = \{\lambda^{(1)}b^{(2)}+\lambda^{(2)}b^{(1)}\}z_1= vz_1/2,$$
 $$\Lambda={v-3 \choose 2}.$$
 The constructed design then has parameters $3-(2v,5,\frac{3}{4}(v-2)(v-4)).$ 
 Since $a_1=(v-2)/6$, we have that $v\equiv 2 \bmod 6.$ 
 Thus we have shown the following.

 \begin{thm}\label{ApplicationIa} 
  There is a simple
  $$ 3-(2v,5, \frac{3}{4}(v-2)(v-4))$$
  design for any integer $v\equiv 2 \bmod 6.$
   
 \end{thm}

\bigskip
 We can construct another family of 3-designs with moderate value for
 $\Theta$. Let $v=2^f+1$ with odd $f$.

\begin{itemize}
 \item
 $D_1$ is the union of $a_1=2^f-1$ copies of
 the complete $2-(2^f+1,2,1)$ design $C_1$. So, $D_1$ is (1,2)-resolvable
 with $\lambda^{(1)}=0$, $\lambda_2^{(1)}=a_1=2^f-1$, $u_1=2$,
 $b^{(1)}=2^f+1$ and $w_1= a_1t_1$ with $t_1=2^{f-1}.$

 \item
 $D_2$ is the complete $3-(2^f+1,3,1)$ design. Since $f$ is odd, we have
 $2^f+1 \equiv 0 \bmod 3.$ So, $D_2$ is $(1,1)$-resolvable.
 For $D_2$ we have $\lambda^{(2)}=1$, 
 $\lambda_2^{(2)}=2^f-1$, $u_2=1$, $b^{(2)}=(2^f+1)/3$ and
 $w_2=2^{f-1}(2^f-1).$  

 \item $D$ is a $3-(2^f+1, 5, 10(2^f-2))$ design, which is obtained
 from the $4-(2^f+1, 5, 20)$ design \cite{Bier-TvT94} with $\gcd(f,6)=1.$
 Thus $\Lambda=10(2^f-2).$

\end{itemize}

 Now
 $$\Theta =\{ \lambda_2^{(1)}u_2+ \lambda_2^{(2)}u_1\}z_1= 3(2^f-1)z_1,$$
 $$\Delta = \{\lambda^{(1)}b^{(2)}+ \lambda^{(2)}b^{(1)}\}z_1= (2^f+1)z_1,$$
 $$\Lambda=10(2^f-2).$$
 Hence  
 $$\Lambda = \Theta - \Delta $$
 if $z_1=5.$ 
 The constructed design has parameters 
 $3-(2(2^f+1),5,15(2^f-1)).$ We have the following.

 \begin{thm}\label{ApplicationIb}
 There is a simple $3-(2(2^f+1),5,15(2^f-1))$ design for $\gcd(f,6)=1.$

 \end{thm}

\subsubsection{\boldmath${n=2}$} \par

\medskip

  We construct a family of simple 3-designs with $k=7$ by using
  Construction I with $n=2.$

 \medskip
 Let $v$ be an integer such that $v\equiv 0 \bmod 4$, $\gcd(v,3)=1$
 and $\gcd(v,5)=1.$
 
\begin{itemize}
 \item
 $D_1$ is the union of $a_1={v-2 \choose 3}/20$ copies of
 the complete $2-(v,2,1)$ design $C_1$. So, $D_1$ is (1,1)-resolvable.
 Here we have $\lambda^{(1)}=0$, $\lambda_2^{(1)}=a_1$, $u_1=1$ and
 $b^{(1)}=v/2$ and $w_1= a_1t_1$ with $t_1=(v-1).$

 \item
 $D_3$ is the complete $3-(v,5, {v-3 \choose 2})$ design, 
 which is $(1,5)$-resolvable.
 For $D_3$ we have $\lambda^{(3)}={v-3 \choose 2}$, 
 $\lambda_2^{(3)}={v-2 \choose 3}$, $u_3=5$, $b^{(3)}=v$ and
 $w_3={v-1 \choose 4}/5.$

 \item
 $D_2$ is the union of $a_2=(v-3)$ copies of the 
 complete $3-(v,3,1)$ design $C_2$. So, $D_2$ is $(1,3)$-resolvable.
 For $D_2$ we have $\lambda^{(2)}=v-3$, 
 $\lambda_2^{(2)}=(v-2)(v-3)$, $u_2=3$, $b^{(2)}=v$ and
 $w_2=a_2t_2$ with $t_2={v-1 \choose 2}/3.$

 \item
 $D_4$ is the complete $3-(v,4, v-3)$ design, which is $(1,1)$-resolvable.
 For $D_4$ we have $\lambda^{(4)}=v-3$, 
 $\lambda_2^{(4)}={v-2 \choose 2}$, $u_4=1$, $b^{(4)}=v/4$ and
 $w_4={v-1 \choose 3}.$  

\end{itemize} 

 We have
\begin{eqnarray*}
 \Theta & = & (\lambda_2^{(1)}u_3+\lambda_2^{(3)}u_1)z_1 +
              (\lambda_2^{(2)}u_4+\lambda_2^{(4)}u_2)z_2 \\
        & = & \frac{5}{4}{v-2 \choose 3}z_1 + 5{v-2 \choose 2}z_2
\end{eqnarray*}

\begin{eqnarray*}
 \Delta & = & (\lambda^{(1)}b^{(3)}+\lambda^{(3)}b^{(1)})z_1 +
              (\lambda^{(2)}b^{(4)}+\lambda^{(4)}b^{(2)})z_2 \\
        & = & \frac{1}{4}v(v-3)(v-4)z_1 + \frac{5}{4}v(v-3)z_2
\end{eqnarray*}

Construction I will yield a simple $3-(2v,7,\Theta)$ design,
when there exist values for $z_1$ and $z_2$ such that
$\Theta -\Delta=0.$  
 
Set  $$ \Theta -\Delta := -Az_1 + Bz_2.$$
Then we have
     $$ A= \frac{1}{24}(v-3)(v-4)(v+10) $$
and
    $$ B=\frac{5}{4}(v-3)(v-4). $$  

It follows that $\Theta -\Delta =0$
if we have $Az_1=Bz_2$, which reduces to the equation
  $$(v+10)z_1 =30 z_2, $$
where $z_1 \leq t_1$ and $z_2 \leq t_2$, i.e.
$z_1 \leq v-1$ and $z_2 \leq (v-1)(v-2)/6.$
It is clear that $z_1=30m$ and $z_2=(v+10)m$ for integer $m \leq (v-1)/30$
are solutions to the equation.
From $z_1=30m$ and $z_2=(v+10)m$ we obtain
$$\Theta =\frac{35}{4}v(v-2)(v-3)m.$$
Recall that $v\equiv 0 \bmod 4$, $v\equiv 1,2 \bmod 3$, and $\gcd(5,v)=1$. 
Moreover, since $a_1={v-2 \choose 3}/20$ must be an integer, we have
$v \equiv 2,3,4 \bmod 5$. Now the congruence system
$v\equiv 0 \bmod 4$, $v\equiv 1,2 \bmod 3$, $v \equiv 2,3,4 \bmod 5$ 
has  $v\equiv 4,8,28,32,44, 52 \bmod 60$ as solutions.
Thus we have proven the following.

\begin{thm}\label{ApplicationIc}
There is a simple $3-(2v,7, \frac{35}{4}v(v-2)(v-3)m)$ design
for any integer $v\equiv 4,8,28,32,44, 52 \bmod 60$ (with $v \geq 32$) 
and  any integer $m \leq (v-1)/30.$
\end{thm}

\bigskip
\subsection{Applications of Construction II}

\subsubsection{\boldmath${n=1}$}

\medskip

 Here is the first example.

\medskip
  Let $f > 3$ be an odd integer such that $\gcd(f,3)=1$.

\begin{itemize}
 \item
 $D_1$ is the complete $3-(2^f+1,3,1)$ design. $D_1$ is $(1,1)$-resolvable.
 For $D_1$ we have $\lambda^{(1)}=1$, 
 $\lambda_2^{(1)}=2^f-1$, $u_1=1$, $b^{(1)}=(2^f+1)/3$ and
 $w_1=2^{f-1}(2^f-1).$  

 \item $D$ is a $3-(2^f+1, 6, \Lambda)$ design, which is obtained
 from the $4-(2^f+1, 6, \lambda)$ design \cite{bier2001} 
 with $\gcd(f,6)=1$, where $\lambda \in \{10,60,70,90,100,150,160\}.$
 Thus $\Lambda=\lambda(2^f-2)/3.$

\end{itemize}

 Now from Theorem \ref{Construction2} we have
 $\Theta^* = \lambda_2^{(1)}u_1z_1,$
 $\Delta^* = \lambda^{(1)}b^{(1)}z_1.$ So,
 $\Theta^* - \Delta^* = \frac{2}{3}(2^f-2)z_1.$
 Thus
 $\Lambda = \Theta^* - \Delta^*$
 if $z_1=\lambda/2.$ 
 The constructed design has parameters 
 $3-(2(2^f+1),6, \Theta^*)$ with $\Theta^*=(2^f-1)z_1= (2^f-1)\lambda/2.$

 We have the following.

 \begin{thm}\label{ApplicationIIa}
 There exists a simple $3-(2(2^f+1),6,(2^f-1)m)$ design for
 $m \in \{5,30,35,45,50,75,80\}$ and $\gcd(f,6)=1.$
 \end{thm}

 We consider another example of general form. 
 Let $v, k$ be integers with $ v > k\geq 3$ and $\gcd(v,k)=1$.

 \begin{itemize}
 \item
 $D_1$ is the complete design $3-(v,k, {v-3 \choose k-3}).$ So,
 $\lambda^{(1)} = {v-3 \choose k-3}$, 
 $\lambda^{(1)}_2 = {v-2 \choose k-2}$, $u_1=k$, $b^{(1)}=v$, and
 $w_1={v-1 \choose k-1}/k.$ 

 \item
 $D$ is a $3-(v,2k, \Lambda)$ design.

 \end{itemize}
 We have
 $\Theta^* = \lambda_2^{(1)}u_1z_1,$
 $\Delta^* = \lambda^{(1)}b^{(1)}z_1.$ 
 Construction II yields a simple $3-(2v,2k, \Theta^*)$ design, when it holds
 $$\Theta^* - \Delta^* =
     (\lambda_2^{(1)}u_1-\lambda^{(1)}b^{(1)})z_1=\Lambda,$$
  or
  $$2{v-3 \choose k-2}z_1 =\Lambda,$$ 
  with $z_1 \leq {v-1 \choose k-1}/k.$ In this case 
  we have 
   $$\Theta^*=\frac{k(v-2)}{2(v-k)}\Lambda. $$ 
 We record the result obtained above.

 \begin{thm}\label{ApplicationII_general}
  Let $v  > k \geq 3$ be integers with $\gcd(v,k)=1.$ Assume that
  there exists a simple $3-(v,2k, \Lambda)$ design such that
  $m=\Lambda/2{v-3\choose k-2}$ is an integer and
  $m \leq {v-1 \choose k-1}/k.$ Then there exists a simple
  $3-(2v,2k, \frac{k(v-2)}{2(v-k)}\Lambda)$ design.
 \end{thm} 

  We will illustrate some
  explicit families for 3-designs from Theorem \ref{ApplicationII_general}
  by taking the
  $3-(v,2k, \Lambda)$ design $D$ to be the complete 
  $3-(v,2k, {v-3 \choose 2k-3})$ design.
 
  \medskip 
  \begin{itemize}
   \item {\boldmath${k=3}.$} 
    $D$ is the  $3-(v,6, {v-3 \choose 3})$ design.
    There exists a simple $3-(2v, 6, \Theta^*)$ design with
    $\Theta^*=\frac{3(v-2)}{2(v-3)}{v-3\choose 3}$,
    if $m={v-3 \choose 2k-3}/2{v-3\choose k-2}=(v-4)(v-5)/12$ is an
    integer. This condition is equivalent to $v\equiv 1,2 \bmod 3$ and 
     $v\equiv 0,1 \bmod 4$. Hence $v\equiv 1,4,5,8 \bmod 12.$
  
    \medskip
   \item {\boldmath${k=4}.$} 
    $D$ is the $3-(v,8, {v-3 \choose 5})$ design.
    There exists a simple $3-(2v, 8, \Theta^*)$ design with
    $\Theta^*=\frac{4(v-2)}{2(v-4)}{v-3\choose 5}$, if
     $m={v-3 \choose 5}/2{v-3\choose 2}=(v-5)(v-6)(v-7)/2.3.4.5$
    is an integer. This condition is equivalent to $v\equiv 1,3 \bmod 4$ and
    $v \equiv 0,1,2 \bmod 5.$ Hence
     $v \equiv 1,5,7,11,15,17 \bmod 20.$

   \medskip
   \item {\boldmath${k=5}.$}
    $D$ is the $3-(v,10, {v-3 \choose 7})$ design.
    There is a simple $3-(2v, 10, \Theta^*)$ design with
    $\Theta^*=\frac{5(v-2)}{2(v-4)}{v-3\choose 7}$, if
     $m={v-3 \choose 7}/2{v-3\choose 2}=(v-6)(v-7)(v-8)(v-9)/16.3.5.7$
    is an integer. This condition is equivalent to $\gcd(v,5)=1$,
    $v\equiv 0,1,6,7 \bmod 8$ and $v \equiv 0,1,2,6 \bmod 7.$

  \end{itemize} 
 In summary, we have the following corollary of Theorem 
 \ref{ApplicationII_general}.

\begin{cor}\label{ApplicationIIb}
 The following hold.

\begin{itemize}
\item[(i)]
There is a simple 
$3-(2v,6, \; \frac{3(v-2)}{2(v-3)}{v-3 \choose 3})$  
design for $v \equiv 1,4,5,8 \; \;  \bmod 12.$

\item[(ii)]
There is a simple 
$3-(2v,8, \; \frac{4(v-2)}{2(v-4)}{v-3 \choose 5})$ 
design for $v \equiv 1,5,7,11,15,17  \; \; \bmod 20.$

\item[(iii)]
There is a simple 
 $3-(2v,10, \; \frac{5(v-2)}{2(v-5)}{v-3 \choose 7})$ 
design for $v \equiv 0,1,2,6 \; \bmod 7,$  $v \equiv 0,1,6, 7 \; \bmod 8,$
 and $\gcd (v,5)=1.$

\end{itemize}

\end{cor}

\subsubsection{\boldmath${n=2}$}

\medskip
 Let $v, k$ be integers such that $ v > 2k$, $k \geq 3$, 
 $\gcd(v,2k)=1$  and $\gcd(v,k+1)=1$.
 
 \begin{itemize}
 \item
 $D_1$ is a union of $a_1=\frac{1}{k(2k-1)}{v-2 \choose 2k-2}$ 
 copies of the complete $2-(v,2, 1)$ design $C_1.$ Since $\gcd (v,2)=1$,
 $C_1$ is $(1,2)$-resolvable and has $t_1=(v-1)/2$ resolution
 classes. For $D_1$ we have 
 $\lambda^{(1)} = 0$, 
 $\lambda^{(1)}_2 = \frac{1}{k(2k-1)}{v-2 \choose 2k-2}$, 
 $u_1=2$, $b^{(1)}=v$, $w_1=a_1t_1.$

\item
 $D_3$ is the complete $2-(v,2k, {v-3 \choose 2k-3})$ design 
 which is $(1,2k)$-resolvable. For $D_3$ we have
 $\lambda^{(3)} = {v-3 \choose 2k-3}$,                     
 $\lambda^{(3)}_2 = {v-2 \choose 2k-2}$, $u_3=2k$, $b^{(3)}=v$, 
 $w_3=\frac{1}{2k}{v-1 \choose 2k-1}$.  

\item
 $D_2$ is the complete $2-(v,k+1, {v-3 \choose k-2})$ design
 which is $(1,k+1)$-resolvable. For $D_2$ we have
 $\lambda^{(2)} = {v-3 \choose k-2}$,
 $\lambda^{(2)}_2 = {v-2 \choose k-1}$, $u_2=k+1$, $b^{(2)}=v$, 
 $w_2=\frac{1}{k+1}{v-1 \choose k}$.

\end{itemize}

 We have
\begin{eqnarray*}
 \Theta^* & = & (\lambda_2^{(1)}u_3+\lambda_2^{(3)}u_1)z_1 +
              \lambda_2^{(2)}u_2z_2 \\
        & = & \frac{4k}{2k-1}{v-2 \choose 2k-2}z_1 + 
               (k+1){v-2 \choose k-1}z_2, \\ \\ \\
 \Delta^* & = & (\lambda^{(1)}b^{(3)}+\lambda^{(3)}b^{(1)})z_1 +
              \lambda^{(2)}b^{(2)}z_2 \\
        & = & v{v-3 \choose 2k-3}z_1 + v{v-3 \choose k-2}z_2
\end{eqnarray*}

We then obtain a simple $3-(2v,2(k+1),\Theta^*)$ design, if
there exist positive integers $z_1$ and $z_2$ 
with $z_1 \leq t_1$ and $z_2 \leq w_2$ for which
$\Theta^* -\Delta^* = 0.$  
 
Set  $$ \Theta^* -\Delta^* := -Az_1 + Bz_2.$$
Then we have
\begin{eqnarray*}
 -A &=& \frac{4k}{2k-1}{v-2 \choose 2k-2}- v{v-3 \choose 2k-3} \\
   &=& -{v-3\choose 2k-3}\alpha
\end{eqnarray*}
with $\alpha =[v(4k^2-10k+2)+8k]/(2k-2)(2k-1),$
\begin{eqnarray*}
   B &=& (k+1){v-2 \choose k-1}-v{v-3 \choose k-2} \\  
     &=& 2{v-3\choose k-2}(v-k-1)/(k-1).
\end{eqnarray*}

Hence, if $\Theta^* -\Delta^* =0$, we have  $Az_1=Bz_2$.
In particular, if $A/B$ is an integer, 
then for any integer $1\leq z_1 \leq t_1$
such that $z_2 =z_1A/B \leq w_2$, we obtain a simple
$3-(2k, 2(k+1), \Theta^*)$ design. 

Here we record this result.

\begin{thm}\label{ApplicationII_general2}
 Let $v, k$ be integers such that $v > 2k$, $k \geq 3$, 
 $\gcd(v,2k)=1$ and  $\gcd(v,k+1)=1.$ 
 Define
$A={v-3\choose 2k-3}\frac{v(4k^2-10k+2)+8k}{(2k-2)(2k-1)}$ and
$B=2{v-3\choose k-2}\frac{v-k-1}{k-1}$. If $A/B$ is an integer, then for
any integer $1\leq z_1 \leq (v-1)/2$ such that $z_2 =z_1A/B \leq
\frac{1}{k+1}{v-1\choose k}$, there exists a simple
$3-(2v,2(k+1), \Theta^*)$ design with
$$\Theta^*={v-2\choose 2k-2}\frac{4k}{2k-1}z_1+{v-2\choose k-1}(k+1)z_2.$$
\end{thm}

\bigskip
We illustrate two special cases with $k=3$ and $k=4$ 
of Theorem \ref{ApplicationII_general2}.

\medskip

\begin{itemize}
\item {\boldmath${k=3}.$}
 
We then have  $A/B=\frac{(v-5)(v-3)}{3.5}$. 
 The conditions that $\gcd(v,6)=\gcd(v,4)=1$ and $A/B$ is an integer 
 are equivalent to $v\equiv 2 \bmod 3$, $v\equiv 1,3 \bmod 4$ and
 $v \equiv 0,2 \bmod 5$. Thus we have
 $v \equiv 5,17,35,47 \bmod 60$. Note that $z_2=z_1A/B.$
 In this case we have a $3-(2v, 8, \Theta^*)$ with
 \begin{eqnarray*} 
 \Theta^*   
  &=& {v-2\choose 4}\frac{12}{5}z_1+{v-2\choose 2}4z_2 \\
  &=& \frac{7}{30}v(v-2)(v-3)(v-5)z_1,
 \end{eqnarray*} 
 where $ 1 \leq z_1 \leq (v-1)/2.$

\medskip
\item {\boldmath${k=4}.$}

 We obtain  $A/B=(v-6)(v-7)(13v+16)/8.3.5.7$ . The requirement that
 $\gcd(v,2k)=\gcd(v,8)=1$, $\gcd(v,k+1)=\gcd(v,5)=1$ and $A/B$ is
 an integer, reduces  to  $v\equiv 7 \bmod 8$, $v \equiv 1,2,3 \bmod 5$
 and $v \equiv 0,2,6 \bmod 7$. Hence
 $v\equiv 7,23,63, 111, 167, 191, 223, 231, 247 \bmod 280.$
  And we have a simple $3-(2v, 10, \Theta^*)$ design with
 \begin{eqnarray*} 
 \Theta^*   
  &=& {v-2\choose 6}\frac{16}{7}z_1+{v-2\choose 3}5z_2 \\
  &=& 81v{v-2\choose 6}z_1/7(v-5).
 \end{eqnarray*} 
 
 \end{itemize}
In summary, we have proven the following.

\begin{cor}\label{ApplII}
 The following hold.
 \begin{itemize}
\item[(i)]
 There is a simple 
 $3-(2v,8, \frac{7}{30}v(v-2)(v-3)(v-5)m)$ design
 for any positive integers  $v \equiv 5,17,35,47 \bmod 60$ and
 $m \leq (v-1)/2.$
\item[(ii)]
There is a simple
$3-(2v, 10,81v{v-2\choose 6}m/7(v-5))$ design
for any positive integers
 $v\equiv 7,23,63, 111, 167, 191, 223, 231, 247 \bmod 280$
 and $ m \leq (v-1)/2.$
\end{itemize}

\end{cor}

 \bigskip
\subsection{\boldmath${(1, \sigma)}$-resolvability of the constructed designs}

\medskip
In this section, we discuss the question of $(1, \sigma)$-resolvability
of the designs obtained by Constructions I and II. In particular, 
we will consider the cases $\Theta-\Delta=0$ and $\Theta^*-\Delta^*=0$,
i.e. the cases where a $3-(v,k,\Lambda)$ design  $D$ is not used in
the construction.

\bigskip
We make use of the following observation. 

\begin{itemize}

\item
Let $(D_h, D_{n+h})$ be a pair of designs
in Constructions I or II such that $k_h \neq k_{n+h}$. 
For given  $(i,j)$ the blocks constructed from the resolution classes 
$(\pi_i^{(h)},{\tilde \pi}_j^{(n+h)})$ and 
$({\tilde \pi}_i^{(h)}, \pi_j^{(n+h)})$ will be denoted by
${\mathcal B}_{h,n+h}^{(i,j)}.$ Thus

$${\mathcal B}_{h,n+h}^{(i,j)}=\{ A\cup \tilde{B},\; \tilde{A}\cup B \; /
        A \in \pi_i^{(h)}, \tilde{A} \in \tilde{\pi}_i^{(h)}, \; 
        B \in \pi_j^{(n+h)}, \tilde{B} \in {\tilde \pi}_j^{(n+h)} \}. $$
Recall that 
      $\varepsilon_h \leq d^{(h)}(\pi_i^{(h)}, \pi_j^{(h)})\leq s_h.$
It follows that each point $x \in X$ or ${\tilde x} \in {\tilde X}$ appears
in 
$$\sigma^{(i)}:= u_hb^{(n+h)}+u_{n+h}b^{(h)}$$
 blocks of ${\mathcal B}_{h,n+h}^{(i,j)}$.
Note that $|{\mathcal B}_{h,n+h}^{(i,j)}|=2b^{(h)}b^{(n+h)}.$

\medskip
\item
For the blocks of type IV in Construction II we have $D_n= D_{2n}$
i.e. $k_n=k_{2n}.$ 
Let ${\mathcal B}_{n,n}^{(i,j)}$ denote the set of blocks
constructed from resolution classes of $D_n$ and 
${\tilde D}_n$ corresponding to the pair $(i,j)$. Then we have 

$${\mathcal B}_{n,n}^{(i,j)}=\{ A\cup \tilde{B}\;  /
        A \in \pi_i^{(n)}, \; \tilde{B} \in {\tilde \pi}_j^{(n)} \}. $$
We have $|{\mathcal B}_{n,n}^{(i,j)}|=b^{(n)}b^{(n)}$ 
and each point $x \in X$ or 
${\tilde x} \in {\tilde X}$ appears in 
            $$\sigma^{(n)}:= u_nb^{(n)}$$
 blocks of ${\mathcal B}_{n,n}^{(i,j)}.$

\end{itemize}
Let $m_1, \ldots, m_n$ be positive integers such that
$$ m_1\sigma^{(1)}= \cdots = m_n\sigma^{(n)} := \sigma .$$
 
Observe that the blocks constructed by each
pair $(D_h, D_{n+h})$ is a union of $z_hw_h$ subsets  
${\mathcal B}_{h,n+h}^{(i,j)}$ of equal size.
Now assume that $m_h | z_hw_h$ for all $h=1, \ldots, n.$
This is equivalent to say that the blocks constructed by the pair
$(D_h, D_{n+h})$ can be partitioned into $z_hw_h/m_h$ disjoint
$1-(2v,k_h+k_{n+h}, \sigma)= 1-(2v,k, \sigma)$ designs.
It is then clear that the constructed design is
$(1, \sigma)$-resolvable. 

\medskip
In summary, by using the notation above we have the following result.

\begin{prop}\label{resolution}
 Let $\mathcal D$ be a $3-(2v,k, \Theta)$ (resp.
  $3-(2v,k, \Theta)^*$) design obtained by Construction I
 (resp. Construction II) for which $\Theta-\Delta =0$
 (resp. $\Theta^*-\Delta^* =0$). 
 Assume that there exist positive integers $m_1, \ldots, m_n$ with
 $m_h | z_hw_h$, for $h=1, \ldots, n$, such that
 $m_1\sigma^{(1)} = \cdots = m_n\sigma^{(n)} := \sigma.$ 
 Then the constructed design $\mathcal D$ is $(1, \sigma)$-resolvable. 
\end{prop}

\bigskip
 In the rest of this section we consider the $(1,\sigma)$-resolvability
 of some families of 3-designs constructed above.
 
 \begin{itemize}

 \item 
 We begin with the
 simple $3-(2v,7, \frac{35}{4}v(v-2)(v-3)m)$ design $\mathcal D$ 
 in Theorem \ref{ApplicationIc}, where
 $v\equiv 4,8,28,32,44, 52 \bmod 60$ (with $v \geq 32$) 
 and integer $m \leq (v-1)/30$.
 The design $\mathcal D$ is obtained by Construction I with $n=2$ 
 and $\Theta -\Delta=0.$ From the parameters of the ingredients
 (see the proof of Theorem \ref{ApplicationIc}) we have

  $\sigma^{(1)}= u_1b^{(3)}+u_3b^{(1)}=v+5v/2=7v/2,$

  $\sigma^{(2)}= u_2b^{(4)}+u_4b^{(2)}=3v/4+ v=7v/4.$ 
  
  Choose $m_1=1$ and $m_2=2$. Then we have 
  $\sigma=\sigma^{(1)}=2 \sigma^{(2)}.$
  Now the condition of Proposition \ref{resolution} reduces to
  $m_2|z_2w_2$, i.e. $2|(v+10)mw_2$, which is always satisfied
  since $v$ is even. Hence $\mathcal D$ is $(1,7v/2)$-resolvable.

 \medskip
 \item
  Consider the designs in Corollary \ref{ApplII} obtained
  by Construction II with $n=2$ and $\Theta^*-\Delta^*=0.$
 
 \begin{itemize}
\item[(i)]  Let $\mathcal D$ be a simple 
 $3-(2v,8, \frac{7}{30}v(v-2)(v-3)(v-5)m)$ design from 
 Corollary \ref{ApplII}, where $v \equiv 5,17,35,47 \bmod 60$ and
 $m \leq (v-1)/2.$  Here we have
  
 $\sigma^{(1)}= u_1b^{(3)}+u_3b^{(1)} =2v+6v = 8v,$

 $\sigma^{(2)}= u_2b^{(2)}=4v.$
 
 Take $m_1=1$ and $m_2=2$, then $\sigma=\sigma^{(1)}=2\sigma^{(2)}=8v.$
 The condition is $m_2|z_2w_2$, i.e. $2|z_2w_2$, where $z_2=z_1A/B$ with 
 $A/B=\frac{(v-5)(v-3)}{3.5}$. Since $v$ is odd, so $A/B$ is even. 
 Thus $2|z_2w_2$. Hence $\mathcal D$ is $(1,8v)$-resolvable.
  
\medskip
\item[(ii)] Similarly,  let $\mathcal D$ be  a simple
$3-(2v, 10,81v{v-2\choose 6}m/7(v-5))$ design from Corollary \ref{ApplII},
with  $v \equiv 7 \bmod 8$, $v\equiv 1,2,3 \bmod 5$,
 $v \equiv 0,2,6 \bmod 7$ and $ m \leq (v-1)/2.$
  We have

$\sigma^{(1)}= u_1b^{(3)}+u_3b^{(1)} =2v+2kv= 10v,$

 $\sigma^{(2)}= u_2b^{(2)}=(k+1)v=5v.$

Take $m_1=1$ and $m_2=2$, then $\sigma=\sigma^{(1)}=2\sigma^{(2)}=10v.$
 The condition is $m_2|z_2w_2$, i.e. $2|z_2w_2$, where $z_2=z_1A/B$ with
 $A/B=\frac{(v-6)(v-7)(13v+16)}{8.3.5.7}$. Thus, if either
 $z_1(=m)$ is even or $A/B$ is even, then the condition
 $2|z_2w_2$ is satisfied. Hence the design $\mathcal D$ is 
 $(1,10v)$-resolvable.
 Note that $A/B$ being an  even integer is equivalent to
 $16|(v-7)$ or $v\equiv 7 \bmod 16$, 
 $v\equiv 1,2,3 \bmod 5$ and $v \equiv 0,2,6 \bmod 7$. 
 
\end{itemize}

\end{itemize}

 We have proven the following.

\begin{prop}\label{ResolvableDesigns}
 The following hold.
\begin{itemize}
\item[(i)]
 The $3-(2v,7, \frac{35}{4}v(v-2)(v-3)m)$ design $\mathcal D$ 
 in Theorem \ref{ApplicationIc} is $(1,7v/2)$-resolvable
 for $v\equiv 4,8,28,32,44, 52 \bmod 60$ (with $v \geq 32$) 
 and integer $m \leq (v-1)/30.$
 
\smallskip

 \item[(ii)]
 The $3-(2v,8, \frac{7}{30}v(v-2)(v-3)(v-5)m)$ design $\mathcal D$ from 
 Corollary \ref{ApplII} is $(1,8v)$-resolvable 
 for $v \equiv 5,17,35,47 \bmod 60$ and $m \leq (v-1)/2.$

\smallskip
 \item[(iii)]
 The $3-(2v, 10,81v{v-2\choose 6}m/7(v-5))$ design $\mathcal D$ 
 from Corollary \ref{ApplII} for
 $v\equiv 7,23,63, 111, 167, 191, 223, 231, 247 \bmod 280$
 and $ m \leq (v-1)/2$ is $(1,10v)$-resolvable, if
 either $m$ even or $16|(v-7).$
 
\end{itemize}

\end{prop}

It is an open question whether Constructions I and II provide 
a $(2,\sigma)$-resolvable 3-design. 

 \medskip
  
  Finally, we include a table listing the 
  simple 3-designs constructed in the paper.

 \medskip

\begin{table}[ht!]
 \caption{\small{ 
      Families of simple 3-designs constructed using Theorems
   \ref{Construction1}, \ref{Construction2}} }

\smallskip
{\small
\begin{tabular}{llll} \hline
 No.    &     Constructed design   & Condition  &  Comment  \\ \hline
    &  &  &  \\
 1   &  $ 3-(2v,5, \frac{3}{4}(v-2)(v-4))$ & $v\equiv 2 \bmod 6$ &
                                      Thm. \ref{ApplicationIa} \\
     &  &  &  \\
 2   & $3-(2(2^f+1),5,15(2^f-1))$  &$\gcd(f,6)=1$ &
                                      Thm. \ref{ApplicationIb} \\
     &  &  &  \\
 3   & $3-(2v,7, \frac{35}{4}v(v-2)(v-3)m)$ & 
              $v\equiv 4,8,28,32,44, 52 \bmod 60$ &
                                    Thm. \ref{ApplicationIc} \\
     &     &  $v \geq 32$, $m \leq (v-1)/30$ &     \\     
     &  &  &  \\
 4   & $3-(2(2^f+1),6,(2^f-1)m)$ & $m \in \{5,30,35,45,50,75,80\}$ &
                                    Thm. \ref{ApplicationIIa} \\ 
     &     &      $\gcd(f,6)=1$  &                       \\
     &  &  &  \\
 5   & $3-(2v,6, \; \frac{3(v-2)}{2(v-3)}{v-3 \choose 3})$ &
             $v \equiv 1,4,5,8  \bmod 12$ & Cor. \ref{ApplicationIIb}(i) \\
     &  &  &  \\
 6   & $3-(2v,8, \; \frac{4(v-2)}{2(v-4)}{v-3 \choose 5})$ &
     $v \equiv 1,5,7,11,15,17  \bmod 20$ & Cor. \ref{ApplicationIIb}(ii) \\
     &  &  &  \\
 7   &  $3-(2v,10, \; \frac{5(v-2)}{2(v-5)}{v-3 \choose 7})$ &
        $v \equiv 0,1,2,6 \bmod 7,$  & Cor. \ref{ApplicationIIb}(iii) \\ 
     &  & $v \equiv 0,1,6, 7  \bmod 8,$ &  \\
     &  &    $\gcd (v,5)=1$  &  \\
     &  &  &  \\
 8   & $3-(2v,8, \frac{7}{30}v(v-2)(v-3)(v-5)m)$ &
     $v \equiv 5,17,35,47 \bmod 60$, & Cor. \ref{ApplII}(i) \\
     &  & $m \leq (v-1)/2$ & \\
     &  &  & \\
 9   & $3-(2v, 10,81v{v-2\choose 6}m/7(v-5))$ &
      $v\equiv 7,23,63, 111, 167,$ &
                                       Cor. \ref{ApplII}(ii) \\  
     &   &  $191, 223, 231, 247 \bmod 280$, &    \\   
     &   &  $m \leq (v-1)/2$   &  \\ 
     & & &  \\ \hline
\end{tabular}
}

\end{table}

\end{document}